\documentclass{emsprocart}

 \usepackage{color}
 
\definecolor{red}{rgb}{1,0,0} 

\contact[k.vogtmann@warwick.ac.uk]{Karen Vogtmann, Mathematics Institute, University of Warwick, Zeeman Building, Coventry CV4 7AL, United Kingdom}

\newtheorem{theorem}{Theorem}[section]

\newcommand{\Q}{\mathbb{Q}}
\newcommand{\R}{\mathbb{R}}
\newcommand{\Z}{\mathbb{Z}}
\newcommand{\G}{\Gamma}
\newcommand{\Out}{{\rm Out}}
\newcommand{\Outn}{{\rm Out}(F_n)}
\newcommand{\Aut}{{\rm Aut}}
\newcommand{\Autn}{{\rm Aut}(F_n)}

\newcommand{\On}{\mathcal {X}_n}% Notation for Outer space of rank n
\newcommand{\Ons}{\mathcal {X}_{n,s}}  %Space for graphs with s leaves
\newcommand{\Mn}{\mathcal M_n} % Notation for the moduli space of graphs.
\newcommand{\Mns}{\mathcal M_{n,s}} % Notation for the moduli space of graphs with s leaves
\newcommand{\OG}{\mathcal{X}_\Gamma} %Notation for Outer space of a RAAG
\newcommand{\FFn}{\mathcal{FF}_n} %Notation for the free factor complex

\title[The topology and geometry of automorphism groups of free groups ]{The topology and geometry of automorphism groups of free groups }

\author[Karen Vogtmann]{Karen Vogtmann\thanks{The author is grateful to
the Royal Society and the Humboldt Foundation for support during the writing of
this paper.}}

\begin{document}

\begin{abstract}
In the 1970's Stallings showed that one could learn a great deal about free groups and
their automorphisms by viewing the free groups as fundamental groups of graphs and
modeling their automorphisms as homotopy equivalences of graphs. Further impetus for
using graphs to study automorphism groups of free groups came from the introduction
of a space of graphs, now known as Outer space, on which the group
$\Outn$
acts nicely.
The study of Outer space and its
$\Outn$
action continues to give new information
about the structure of
$\Outn$, but has also found surprising connections to many other
groups, spaces and seemingly unrelated topics, from phylogenetic trees to cyclic operads
and modular forms. In this talk I will highlight various ways these ideas are currently
evolving.
\end{abstract}



\maketitle

\section{Introduction}

The most basic examples of infinite finitely-generated groups are free groups, with free abelian groups and fundamental groups of closed surfaces close behind.  The automorphism groups of these groups are large and complicated, perhaps partly due to the fact that the groups themselves have very little structure that an automorphism needs to preserve.   These automorphism groups have been intensely studied both for their intrinsic interest and because of their many ties to other areas of mathematics.  Inner automorphisms are well understood, so  attention often turns to the outer automorphism group $\Out(G)$, i.e. the quotient of $\Aut(G)$ by the inner automorphisms.

The outer automorphism group $\Out(\Z^n)={\rm{GL}}(n,\Z)$ of a free abelian group has classically been studied via its action on the symmetric space $\rm{SL}(n,\R)/\rm{SO}(n)$, and the outer automorphism group $\Out(\pi_1(S))=\rm{Mod}(S)$ of a surface group via its action on Teichm\"uller space $\mathcal T_S$.   These spaces are contractible manifolds with a proper action by $\Out(G).$ Both the geometry of the space and the topology of its quotient by the action yield algebraic information about the group $\Out(G)$.  

  For the free group $F_n$   an analogous space $\On,$  now known as ``Outer space," was introduced by Culler and the author in \cite{CuVo}.  It is a contractible space with a proper action by $\Outn$ but is not a manifold.  Nevertheless it can be endowed with a useful metric structure, and again both the geometry of the space and the topology of the quotient yield information about $\Outn$.   
  
  Outer space has a rich combinatorial structure, based on the fact that its points correspond to  finite graphs with fundamental group $F_n$.  Although Outer space was introduced to study $\Outn,$ the fact that finite graphs are used to parameterize many phenomena in mathematics and science means that the structure of Outer space is connected with such diverse areas as  the study of Lie algebras of derivations, degenerations of algebraic varieties, the computation of Feynman integrals, and the statistics of phylogenetic trees. Several of these connections are indicated briefly in the final section of this article.

This article is meant as a brief introduction to this area of mathematics. The area is very active and evolving rapidly, so it is not possible to give a complete survey.  No proofs are presented, only constructions and statements, but there are many references included for readers who want more details and information.  After a review of the basics I discuss a few of the more recent developments.  These include: 
\begin{enumerate}
\item Variations and generalizations  of  Outer space.  These include Auter space and spaces of graphs with more leaves, deformation spaces,  sphere systems in 3-manifolds, and  spaces of CAT(0) cube complexes.  
\item    New understanding of the cohomology of the quotient $\On/\Outn$ via the combinatorics of the spine and connections to the cohomology of the Lie algebra of symplectic derivations of a free Lie algebra.
\item New understanding of what's at infinity for Outer space, and new ways to use that information.
\end{enumerate}

\section{The Basics}\label{sec:basics}

\subsection{Definition of Outer space} \label{ss:defs} There are number of equivalent definitions of Outer space $\On$, each with its own advantages. Historically the first description was as a space of metric graphs with fundamental group $F_n$. We begin by making this definition precise.

{\bf Points.} By  a {\em metric graph} we mean a finite connected graph with positive real edge lengths, equipped with the path metric.     
We fix a model rose $R_n$ (a graph with one vertex and $n$ petals),  and identify the petals of $R_n$ with the generators of the free group $F_n$.
A point in $\On$ is then a metric graph $G$ together with a homotopy equivalence $g\colon R_n\to G$ called a {\em marking}; the marking serves to identify the fundamental group of $G$ with $F_n$.  Marked graphs $(g,G)$ and $(g',G')$ are considered the same if there is an isometry $f\colon G\to G'$ with $f\circ g$ homotopic to $g'$.  

To get a finite-dimensional space we assume $G$ has no univalent or bivalent vertices; the Euler characteristic then tells us that there are only a finite number of possible combinatorial types of graphs (they have at most $3n-3$ edges).  It is also often convenient to  normalize our objects, which we can do by assuming  that the sum of the edge lengths is equal to $1$ (or, equivalently, consider projective classes of marked metric graphs).  A further constraint which is often convenient is to assume $G$ has no separating edges; this subspace is sometimes called {\em reduced} Outer space.

{\bf Topology.} In order to call $\On$ a ``space" we need to define a topology. Intuitively, a neighborhood of $(g,G)$ is obtained by varying the edge lengths of $G$, but there is a wrinkle:  if a vertex $v$ has valence at least four then nearby points can be obtained by replacing $v$ by a tiny tree and attaching the edges adjacent to $v$ to the tree's leaves. Collapsing the tree again recovers $(g,G)$.  More formally, each marked graph $(g,G)$ determines an open simplex $\sigma(g,G)$, obtained by varying the (positive) edge lengths, keeping the sum equal to $1$.    The simplex $\sigma(g',G')$ is a {\em face} of $\sigma(g,G)$ if $(g',G')$ can be obtained from $(g,G)$ by collapsing some edges to points.  Then $\On$  is the quotient space obtained from the disjoint union of the open simplices $\sigma(g,G)$ by face identifications.  

{\bf Simplicial closure.}  Note that not all faces of the simplex $\sigma(g,G)$ are in $\On$.  For example a vertex of $\sigma(g,G)$ corresponds to a graph with exactly one non-zero edge, but such a graph   does not have fundamental group $F_n$ so does not belong to $\On$.  
If we replace each open simplex $\sigma(g,G)$ by a closed simplex $\overline\sigma(g,G)$ and take the disjoint union modulo face relations as before,  the resulting quotient space {\em is} a simplicial complex $\On^*$ called the {\em simplicial closure} of Outer space.  The points of $\On^*$ which are not in $\On$ are said to be {\em at infinity}

{\bf Action.} The group $\Outn$ acts on $\On$ by changing the marking:   any automorphism $\varphi$ can be realized by a homotopy equivalence $f\colon R_n\to R_n$, and the action of $\varphi$ on $(g,G)$ is then given by  $$(g,G)\varphi=(g\circ f,G).$$
The stabilizer of $(g,G)$ under this action is isomorphic to the (finite) group of isometries of $G$.

{\bf Moduli of graphs.} The quotient $\Mn=\On/\Outn$ is called the {\em moduli space of graphs}.  By a classical result of Hurewicz, the cohomology of $\Mn$ is isomorphic to the algebraically-defined group cohomology of $\Outn,$  when both cohomologies are taken with trivial rational coefficients.  

$\Outn$ contains torsion-free  subgroups of finite index.  Since the action of $\Outn$ on $\On$ has finite stabilizers, any torsion-free subgroup $\Gamma$ acts freely, so the quotient $\On/\Gamma$ is an actual $K(\Gamma,1)$, and the homology with any coefficients of $\On/\Gamma$ is equal to the cohomology of $\Gamma$.  

{\bf The spine of Outer space.}   Since graphs in $\On$ have at most $3n-3$ edges, whose lengths must sum to $1$, the dimension of $\On$ is $3n-4$.  Thus we see immediately that the rational homology of $\Outn$ vanishes in dimensions larger than $3n-4$.  But in fact $\On$ contains an equivariant deformation retract $K_n$ called the {\em spine of Outer space}, which has dimension $2n-4$.  The spine $K_n$ is a subcomplex of the barycentric subdivision of the simplicial closure $\On^*,$ consisting of simplices spanned by vertices which are not at infinity.  In other language, $K_n$ is  the geometric  
realization of the partially ordered set of open simplices $\sigma(g,G)$ in $\On$, where the partial order is given by the face relations.

{\bf Actions on trees.} 
 Given a marked graph $(g,G)$, the metric on $G$ lifts to a metric on the universal cover $\widetilde G$ with the property that the fundamental group $\pi_1(G)$ acts on $\tilde G$ by isometries. The marking $g$ identifies $\pi_1(G)$ with $F_n$, so to each point of  $\On$  we can associate a free isometric action of $F_n$ on a metric tree.  The constraints we have placed on $G$ imply that this action is {\em minimal}, i.e. there are no invariant subtrees.  Conversely, given a minimal action of $F_n$ by isometries on a metric simplicial tree $T$,  the quotient of $T$ by the action is a graph that naturally comes with a marking, determined by the images of paths from an arbitrarily chosen point of $T$ to its  translates  by the generators of $F_n$.   

Any set of isometric actions of a group on metric spaces can be given the {\em equivariant Gromov-Hausdorff topology}.  In this topology, actions of $G$ on spaces $X$ and $Y$ are close if for every $\epsilon>0$, every finite set $A$ of group elements and  every finite set $\{x_1,\ldots,x_m\}$ of points in  $X$, there is a corresponding set of points $\{y_1,\ldots,y_m\}$ in $Y$  so that   $d_X(x_i,ax_j)$ is $\epsilon$-close to  $d_Y(y_i,ay_j)$  for all $a\in A$ and $i,j=1,\ldots,m$.  For $\On$, the quotient topology defined earlier is equivalent to the equivariant Gromov-Hausdorff topology.

{\bf Sphere complexes.} 
  A quite different construction based on embedded $2$-spheres in a certain $3$-manifold, yields another useful description of $\On$. The $3$-manifold $M_n$ in question is the connected sum  of $n$ copies of $S^1\times S^2$.   
By a theorem of Laudenbach \cite{Lau}, the mapping class group  of $M_n$ is an extension of $\Outn$ by a finite $2$-group, where the $2$-group is generated by Dehn twists along embedded $2$-spheres.   Such a Dehn twist does not affect the isotopy class of any embedded $2$-sphere, so the group $\Outn$ acts on the simplicial complex whose vertices are isotopy classes of non-trivial embedded $2$-spheres, and whose $k$-simplices are disjointly embedded sets of $k+1$ (isotopy classes of) $2$-spheres.  This complex is called the  
 {\em sphere complex of $M_n$}, and it turns out to be isomorphic to the simplicial closure $\On^*$.
 To recover $\On$ itself, put barycentric coordinates on the simplices of the sphere complex, so that a point is a sphere system together with non-negative real weights on the spheres, adding up to $1$.   A point is in $\On$ if and only if the spheres with positive weights cut $M_n$ into simply-connected pieces.

\subsection{Contractibility}\label{ss:contract}
The theorem which cements the relationship between $\On$ and $\Outn$ was proved by M. Culler and the author in 1986:  

\begin{theorem}[\cite{CuVo}] $\On$ is  contractible  and the action of  $\Outn$ is proper.  The spine $K_n$ is an equivariant deformation retract of  dimension $2n-3$ with compact quotient.  
\end{theorem}

 There are nice proofs of this theorem from all three points of view, i.e. using the descriptions of $\On$ as a space of marked graphs \cite{CuVo, V_Rep},  a space of actions on trees \cite{Sko, GiLe1}, or a space of weighted  sphere systems \cite{Hat, HaVo1}. The proofs using sphere systems and actions on trees also prove that the simplicial closure $\On^*$ is contractible, and the proof using actions on trees shows that an even larger space of actions (called {\em very small} actions) is contractible.  The set of very small actions can be regarded as a compactification of $\On$; a few more words about this  can be found in Section~\ref{ss:projective}.

\section{Variations and generalizations}\label{sec:variations}

There are many ways to modify the definition of $\On$ to obtain other useful spaces.  In this section we briefly mention a few of the most common. 
 
\subsection{Auter space}\label{ss:auter}
In the graph description of $\On$ the graphs do not come with basepoints, so  markings, isometries and homotopy equivalences  obviously cannot be required to preserve basepoints.  This is the reason that the action is by {\em outer} automorphisms of $F_n$.  To get a space with  an action by $\Autn$ instead of $\Outn$, we simply modify the definition to include basepoints, and insist that all maps preserve them.  Since $\Autn$ acts on this space, Frederic Paulin dubbed it {\em Autre espace}, which has since became anglicized to  {\em Auter space}.  There is a natural map from Auter space to Outer space that  forgets the basepoints.  The fiber of this map over a point $(g,G)$ is the set of  basepointed paths in $G,$ up to homotopy relative to the endpoints, i.e. it is the (contractible) universal cover of $G$; this is one route to proving that Auter space is itself contractible.

Sphere complexes also offer a natural way to define Auter space:  if one ``punctures" the 3-manifold $M_n$ by removing a small ball  the sphere complex of the resulting $3$-manifold is isomorphic to the simplicial closure of Auter space; the puncture corresponds to the basepoint of a graph. The map from Auter space to $\On$ fills in the puncture, so that spheres which were prevented only by the puncture from being isotopic now become isotopic.    The sphere complex proof of contractibility extends without change to sphere systems in the punctured $3$-manifold.  

From the point of view of  tree actions Auter space consists of actions on trees with specified basepoints.   One of the more delicate points in   tree-action proofs that Outer space $\On$ is contractible is showing that one can choose basepoints in a continuous way, so having basepoints as part of the structure actually makes the proof of contractibility considerably simpler.

\subsection{Graphs with leaves}\label{ss:leaves}
There are many situations in which it is desirable to consider graphs with more than one distinguished point.   One way to do this is to allow graphs to have edges called {\em leaves}  which terminate in univalent vertices; the initial vertex of a leaf can then be thought of as ``distinguished point" on the core graph.  Note that this construction allows one to consider graphs where distinguished points ``collide."  The model rose $R_n$ in the graph defiition of $\On$ is replaced by a {\em thorned rose} $R_{n,s}$, which is a graph with $n$ loops at a basepoint and $s$ leaves emanating from the basepoint, labeled by the numbers $1,\ldots,s.$   As before, a point is a marked metric graph, but we  don't assign lengths to the leaves, since sometimes we may want to think of them as distinguished points (length 0 leaves), or sometimes as infinite rays. All maps must send univalent vertices to univalent vertices with the same label.   The space $\Ons$ of marked metric graphs with fundamental group $F_n$ and $s$ labeled leaves is contractible and the group $\Gamma_{n,s}$ of homotopy equivalences of $R_{n,s}$ which fix the univalent vertices  acts properly on $\Ons$.   

The sphere complex point of view is very natural here; all we are doing is allowing several punctures (instead of a single puncture, as we did for Auter space); the definition of the sphere complex and proof of contractibility apply with no significant changes.  

\subsection{Tropical curves}\label{ss:tropical}
Tropical geometers, starting from a different perspective grounded in algebraic geometry are led to define moduli spaces of metric graphs with leaves, which they call {\em tropical curves} (in this theory the ``leaves" are infinite rays, see \cite{Cap}).   
In our terminology this moduli space is the quotient of the space $\Ons$ by the group $\G_{n,s}$ of (homotopy classes of) homotopy equivalences which fix the univalent vertices.    It is desirable in the context of algebraic geometry to compactify moduli spaces by allowing degenerations.  In the case of tropical curves, this means allowing connected subgraphs to collapse to points, retaining only the rank of the collapsed subgraph as a decoration on the resulting vertex.  This translates in our language to taking the quotient of the simplicial closure $\Ons^*$ by the groups $\Gamma_{n,s}$ instead of the quotient of $\Ons$ itself.  

\subsection{Free products, deformation spaces}\label{deformation}
Guirardel and Levitt greatly generalized the definition of $\On$ as a space of actions on trees.  They first defined a space of actions of a free product $G_1\ast\ldots\ast G_k$ on trees, which gave information about the outer automorphism group of this free product in terms of the outer automorphism groups of the factors \cite{GiLe1}.  They then generalized this further to define the notion of a {\em deformation space} of actions of a given group on trees with stabilizers in a fixed class of subgroups \cite{GiLe2}; deformation spaces generalize not only the Outer space of a free product but also  the space of ``JSJ decompositions" of a given group $G$; these are special splittings of $G$ which encode all possible splittings, and are a key tool in geometric group theory.

\subsection{Sphere complexes}\label{ss:spheres}
Hatcher and Wahl generalized the sphere complex definition of $\On$ by considering complexes of isotopy classes of sphere systems in quite general $3$-manifolds $M^3$, in order to study the homology of mapping class groups $\pi_0({\rm Diff}(M^3))$ \cite{HaWa}.  Their main result shows that the inclusions given by repeatedly taking connected sum with a fixed closed $3$-manifold  eventually induce isomorphisms on the $i$-th homology of the associated mapping class groups. This theorem  unifies and generalizes homological stability results about $\Autn$ and $\Outn$, as well as several other series of groups important in low-dimensional topology, such as braid groups and handlebody subgroups of surface mapping class groups.      

\subsection{RAAGs}\label{ss:RAAGs}
The free group $F_n$ is an example of a {\em right-angled Artin group}, usually abbreviated as  RAAG.  Another example is the free abelian group $\Z^n$.   In general a RAAG is a finitely-generated group, defined by the condition that some pairs of its generators commute.  This is usually codified by drawing a graph $\Gamma$ with one vertex for each generator and one edge between each pair of commuting generators; the resulting group is denoted $A_\Gamma$.   In \cite{ChVo} Charney and the author defined an ``Outer space $\OG$ for RAAGs,"  where the model space (the domain of the markings) is a {\em Salvetti complex} $S_\Gamma$;  this is a cube complex with one vertex, an edge (i.e. a loop) for every generator of $A_\Gamma$, and a $k$-cube ($k$-torus) for every set of $k$ mutually commuting generators.  The fundamental group of $S_\Gamma$ is identified with $A_\Gamma$, and any automorophism of $A_\Gamma$ can be realized as a homotopy equivalence of $S_\Gamma$.  

A general point of $\OG$ is a pair $(g,B)$, where $g\colon S_\Gamma \to B$ is a homotopy equivalence and $B$ is a type of locally $CAT(0)$ cube complex  called a  {\em $\Gamma$-complex};  it can be constructed from $S_\Gamma$ by successively sliding various cubes halfway across other cubes.  There are constraints on which slides are allowed, for example they must preserve the fundamental group and the property of being a $CAT(0)$ cube complex.    The restriction of the metric on $B$ to each cube is flat but the sides of the cube are not necessarily orthogonal; some ``cubes" are allowed to be more general parallelepipeds.    The space of all such marked metric objects $(g,B)$ is a hybrid of the space $\On$   of marked graphs and the symmetric space for ${\rm GL_n},$ which can be described as a space  of marked flat $n$-tori.  This outer space is conjectured to be contractible; this was proved for the subspace on which all cubes are orthogonal \cite{ChVo}; in particular if there are no ``twists" in ${\rm Out}(A_\Gamma)$ the space is contractible, and it is known to be contractible in various other cases, such as when the subgroup generated by twists is normal.   The conjecture for all RAAGs is still open as of this writing.   

\section{Algebra from topology and geometry}\label{sec:algebra}

We note here a few ways that having a nice space with a proper action is useful for computing algebraic invariants of $\Outn$.  

\subsection{Topology}\label{ss:topology}
Since the action of $\Outn$ on $\On$ is proper any torsion-free finite-index subgroup $H$ of $\Outn$ acts freely on $\On$. Since $\On$ is contractible, a theorem of Hurewicz implies that the homology of the quotient $\On/H$ can be identified with the algebraically-defined group homology of $\On$.  If one considers homology with trivial rational coefficients then the  homology of the entire group $\Outn$ can be identified with the homology of the quotient $\On/\Outn$; the heuristic reason is that the homology with rational coefficients of a finite group is trivial, so rational cohomology does not  ``see" the finite stabilizers and thinks the action is free.  Thus the topology of the quotient space provides algebraic invariants (such as cohomology) for the group.

\subsection{Geometry}\label{ss:geometry}

The classical Milnor-Svarc Lemma says that a finitely-generated group which acts properly and cocompactly on a simply-connected space $X$ is {\em quasi-isometric} to $X$. This means that if you choose a point $x\in X$, then for any  group elements $g$ and $h$, the distance between   $gx$ and $hx$  is approximately the same as the distance between $g$ and $h$ in the word metric of the group.

The action of $\Outn$ on $\On$ is proper but not cocompact, so $\Outn$ is not quasi-isometric to $\On$.   One way of resolving this is to use the spine $K_n$, which is cocompact,  with its natural simplicial metric.  Thus quasi-isometry invariants for $K_n$ such as the number of ends,  isoperimetric functions in any dimension (including the 2-dimensional Dehn function) or the asymptotic dimension provide quasi-isometry invariants of the group $\Outn$.  
This  was used, for example, to establish upper bounds for isoperimetric functions   in all dimensions  using the sphere complex description of $K_n$  \cite{HaVo2}.  These upper bounds are all exponential in $n$, and in \cite{BrVo, HaMo2} it was shown that for the Dehn function at least there is also an exponential lower bound.  

In recent years a great deal of work has been done on a geometric theory for the entire space $\On$ (as opposed to the spine), where $\On$ is given an asymmetric metric called the {\em Lipschitz metric.}  A few more words about this theory and its consequences   can be found in Section~\ref{sec:Lipschitz} below.  

\subsection{Fixed points}\label{fixed} Any finite subgroup $H$ of $\Outn$ can be realized as the group of automorphisms of a graph with fundamental group $F_n$; this can be rephrased as the statement that  the fixed point set of $H$ is non-empty.  In fact the fixed point set of any finite subgroup is contractible \cite{Whi}, a property shared with the action of mapping class groups on Teichm\"uller spaces and the action of various arithmetic groups on symmetric spaces.  A contractible space with a proper action, where the fixed point sets of finite subgroups are   contractible and fixed point sets of infinite subgroups are empty is called a {\em classifying space for proper actions}.     Such a space is universal among spaces with proper $G$-actions, so plays a central role in many topology-based arguments in group theory.   
 
\section{The topology of moduli spaces of graphs}\label{sec:topology}
 The quotient of $\On$ by $\Outn$ is called the {\em moduli space of graphs} and denoted $\mathcal{MG}_n$. More generally, if we consider graphs with $s$ leaves we write  $\Mns$; this is a rational classifying space for the group $\G_{n,s}$ of homotopy equivalences fixing the univalent vertices of such a graph.  To study topological invariants such as cohomology of $\Mns$ we can use the fact that the spine $K_{n,s}$ of $\Ons$ is a simplicial complex and the stabilizer of a simplex under the action of $\G_{n,s}$  fixes it pointwise.  Thus the simplicial structure of $K_{n,s}$ induces a simplicial cell structure on the quotient space $K_{n,s}/\G_{n,s}$, which we can use to build a chain complex.

\subsection{Cubical structure of the spine}
For simplicity we stick to the case $s=0$ in this paragraph and the next. Recall that the spine $K_{n}$ is the geometric realization of the partially ordered set of open simplices $\sigma(g,G)$.  In other words, there is one vertex for each marked  graph,  where we ignore the metric on $G$.  Vertices $(g,G)$ and $(g',G')$ are joined by an edge if $G'$ can be obtained from $G$ by collapsing some edges to points and $g'$ is homotopic to the composition of $g$ with the collapsing map. Implicit in the statement that $g'$ is a homotopy equivalence is the fact that the set of collapsing edges cannot contain a cycle, i.e. it is a {\em forest} (a disjoint union of trees) so the collapsing map is called a {\em forest collapse}.  A chain of $k$ forest collapses gives a $k$-simplex.   

Given a forest $\Phi$ in $G$ with $k$ edges, we can get a chain of $k$ forest collapses by  by collapsing the edges of $\Phi$ one at a time, in some order. Collapsing in a different order gives a different $k$-simplex,  and all of the $k$-simplices obtained this way fit together to triangulate a $k$-dimensional cube.  Thus we may also think of $K_n$ as a cube complex, with one cube for each triple $(g,G,\Phi)$.  %This gives rise to a chain complex  with many fewer generators, as follows.  %, where $(g,G)$ is a topological marked graph and $\Phi$ is a forest in $G$.    

Each cube $(g,G,\Phi)$ can be oriented by choosing an ordering of the edges of $\Phi$, up to even permutation.  The stabilizer of $(g,G,\Phi)$  is isomorphic to the group of automorphisms of the pair $(G,\Phi)$, and acts linearly on the boundary of the cube.  The quotient by this action is either contractible or is a rational homology sphere, depending on whether there is an orientation-reversing automorphism of the pair $(G,\Phi)$.  Thus the rational homology of $K_n/\Outn$ (and therefore of $\Outn$) can be computed from a chain complex with one generator for each pair $(G,\Phi)$ that has no orientation-reversing automorphisms.    One advantage of this chain complex is that it has many fewer generators than the simplicial chain complex.  Another is that it is very closely related to Kontsevich's graph homology chain complex.  

\subsection{Euler characteristic and rational Euler characteristic}
The cube complex description of $K_n$ and the resulting relatively small chain complex for the quotient by $\Outn$ make it feasible to compute the Euler characteristic of $\Outn$ for small values of $n$ with the help of a computer (see    \cite{HaVo0}).  In fact the Euler characteristic of $\Outn$ has been computed for $n\leq 11$ by Morita, Sakasai and Suzuki (though by other methods) \cite{MSS}.  The values are $1$ or $2$ in ranks less than $9$; the remaining three values are $-21$  (n=9), $-124$ (n=10)  and $-1202$ (n=11), suggesting that for $n$ large $\Outn$ probably has a lot of rational homology, much of it in odd dimensions.  

The {\em rational Euler characteristic} of $\Outn$ is defined as the Euler characteristic of a torsion-free finite-index subgroup divided by the index of the subgroup; this is independent of the choice of finite-index subgroup.  The rational Euler characteristic   is easier to compute and better behaved than the actual Euler characteristic; there is a generating function for it \cite{SmVo}, and the values have been computed  for $n\leq 100$ (though only published for $n\leq 11$).  These values are all negative and seem to grow faster than exponentially.  

The values of the rational Euler characteristics for $n=9, 10, 11$ are approximately $-29, -206, -1691$, which are  tantalizingly close to the actual Euler charactistic values $-21, -124, -1202$.  A natural conjecture, hinted at by Kontsevich and made explicitly by Morita, Sakasai and Suzuki, is that these two Euler characteristics are asymptotically the same.  Heuristically this is because the difference in the two calculations involves the automorphism groups of graphs with fundamental group $F_n$, and one might expect that ``most" such graphs have no automorphisms.

\subsection{Stable homology}
There are natural inclusions from $\Gamma_{n,s}$ into $\G_{n+1,s}$ and $\G_{n,s+1}$, and it is known that the map induced on the $i$-th homology groups by these inclusions is an isomorphism for $n$ sufficiently large with respect to $i$; the same is true for the  map $\G_{n,1}=\Autn\to \G_{n,0}=\Outn$ induced by the quotient map (\cite{HaVo2, HVW}).  Galatius has shown that the stable rational homology is trivial, and in fact the inclusion of the symmetric group into $\Autn$ induces an isomorphism on homology (even with integral coefficients) for $n$ sufficiently large \cite{Gal}.   

The dimension of the spine of $\Mns$ is  $2n-3+s$ for all $n$ and $s$, so it is also clear that the rational homology $H_i$  vanishes for $i>2n-3+s$.  Thus  the large number of rational classes which must exist according to the Euler characteristic calculations live in dimensions   $4n/5 -1 \leq i \leq 2n-3+s$.

\subsection{Unstable classes and assembly maps}
The first non-trivial rational homology class was found by Hatcher and the author in \cite{HaVo0}.  It lies in $H_4(\Aut(F_4))$ and survives under the map to $H_4(\Out(F_4))$.    Later Morita  %, relying on work of Kontsevich \cite{K1,K2},
 found an infinite sequence of cocycles representing potentially nontrivial cohomology classes $\mu_k\in H^{4k}(\Out(F_{2k+2}))$ and showed that the first one is non-zero \cite{Mor}. Morita's construction relies on work of Kontsevich, one of whose {\em graph homology} theories identitified the cohomology of $\Outn$ with the homology of an infinite-dimensional symplectic Lie algebra $\ell_\infty$ \cite{Kon1,Kon2}.  Morita's cocycles are pullbacks of certain elements of the abelianization of $\ell_\infty$, namely those in the image of his {\em trace map}.  

The $\mu_k$ are now called {\em Morita classes}, and Morita has conjectured that they are all non-zero as elements of cohomology. Conant and the author reinterpreted and generalized these classes in \cite{CoVo2}, and showed that $\mu_2$ is also non-zero; Gray later extended this to show $\mu_3$ is non-zero \cite{Gra}. The proofs for $\mu_2$ and $\mu_3$ rely partly on computer calculations which become large extremely fast with $n$, and Morita's conjecture remains a challenging open problem. 

Another non-trivial rational homology class, this time in $H_7(\Aut(F_5))$ was found by Gerlits \cite{Ger}, and also relied on the help of a computer.  This class does not survive the map to $H_7(\Out(F_5))$, and at first did not seem to have any relation to the Morita picture.   This changed in 2011, when Conant, Kassabov and the author introduced a ``hairy" version of Kontsevich's graph homology \cite{CKV}.  They used this to generalize Morita's trace map and find new pieces of the abelianization of $\ell_\infty$.  These new pieces are closely related to modular forms for $\rm{SL}(2,\Z)$, and can be used to construct new cycles; for example, those coming from Eisenstein series can be used to construct cycles in $H_{4n+3}(\Aut(F_{2n+3}))$.  The first of these {\em Eisenstein classes} lies in $H_7(\Aut(F_5))$ and can be identified with Gerlits' class.  The second Eisenstein class, in $H_{11}(Aut(F_7)),$ is also known to be non-trivial.

The unstable homology picture was reformulated, simplified and extended by Conant, Hatcher, Kassabov and the author in \cite{CHKV}.  This paper  avoids the symplectic Lie algebra altogether and introduces a new construction that builds classes in $H_i(\G_{n,s})=H_i(\Mns)$ by ``assembling" classes from the homology of  groups $\Gamma_{m,t}$ associated to a system of subgraphs of lower rank.  The Morita classes are assembled from two rank one classes, and the Eisenstein classes from a rank one class and a rank two class.  As $n$ grows, the number of ways of decomposing a graph of rank $n$ into smaller rank graphs also grows, very fast, so it is plausible that assembling classes could account for the  rapid growth seen in  the Euler characteristic.  

 Until January 2016 this assembly construction accounted for all known non-trivial homology for $\Autn$ and $\Outn$. However, at that time Bartholdi managed to extend the computer calculations to rank $7$ and discovered two new rational homology classes, in $H_8(\Out(F_7))$ and $H_{11}(\Out(F_7))$ \cite{Bar}.  This was unexpected;  note that the Euler characteristic  does not see these classes because one is even-dimensional and one is odd-dimensional.   It is likely that the class in dimension $8$ can be assembled from rank one classes (in a tetrahedral pattern), but the class in dimension $11$ is more surprising:  there is no obvious candidate for a representative in the image of an assembly map and it is the first class found in the virtual cohomological dimension of $\Outn$.  Mapping class groups and   ${\rm{SL}}(n,\Z)$ have no rational homology in their virtual cohomological dimension, so this is also unexpected if one is led by the philosophy that there is a strong analogy between $\Outn$ and these groups.   

\section{The simplicial closure and its quotient}

\subsection{Hyperbolicity of the simplicial closure and related complexes} \label{ss:hyperbolic}
The group $\Outn$ is not Gromov hyperbolic  since it contains free abelian subgroups of large rank $(2n-3)$.  Thus $\Outn$ cannot act properly and cocompactly on any hyperbolic space.  This echoes the situation for mapping class groups $Mod(S)$, but in the case of mapping class groups Masur and Minsky proved that a closely related space with a natural action, the {\em curve complex} is hyperbolic, and they proceeded to use that fact to derive new information about mapping class groups.

There are several candidates for an analog of the curve complex in the setting of $\Outn$. The first is the simplicial closure $\On^*$ of Outer space, which can also be interpreted as the sphere complex of a doubled handlebody, or as the complex of free splittings of a free group.  Using this last characterization Handel and Mosher proved that $\On^*$ is in fact Gromov hyperbolic \cite{HaMo2}; a simpler version of their proof using the language of sphere complexes was later given by Hilion and Horbez \cite{HiHo}.  

The action of $\Outn$ on  $\On^*$ is missing some desirable features enjoyed by the action of  $Mod(S)$ on the curve complex of $S$. For example a mapping class is pseudo-Anosov if and only if it has positive translation length on the curve complex, whereas elements of $\Outn$ with positive translation length in $\On^*$ are not necessarily  {\em fully irreducible}, 
a property often considered analogous to the pseudo-Anosov property for mapping classes.    
There is, however,  a different complex, the {\em free factor complex} $\mathcal{FF}_n$, whose vertices are conjugacy classes of free factors of $F_n$.  This complex too is hyperbolic, as shown by Bestvina and Feighn  \cite{BeFe}, who also showed that an element of $\Outn$  has positive translation length in $\mathcal{FF}_n$ if and only if it is fully irreducible.   
 
It is conjectured that the action of $\Outn$ on $\mathcal{FF}_n$ is {\em acylindrical}, which is a weak analog of proper discontinuity.  Proper discontinuity says that  only finitely many group elements ``almost fix" any single point, while acylindricity  says that the set of group elements which almost fix a far-apart pair of points is finite.   (More formally,   an action is acylindrical if   
given any $\epsilon>0$ there are numbers $R$ and $N$ such that if $x$ and $y$ are distance at least $R$ apart then at most $N$ group elements can displace both $x$ and $y$ by less than $\epsilon$.)  Although it is not known whether the action on $\mathcal{FF}_n$ is acylindrical, it is true that   $\Outn$ is an {\em acylindrically hyperbolic} group; this is proved by constructing a quasi-tree using translates of the axis of a fully irreducible element acting on Outer space (see \cite{BBF} for details).   This property turns out to be strong enough for  many geometric arguments about the group; see, e.g. \cite{Osi} for a survey.   

\subsection{Moduli spaces of graphs with leaves}\label{ss:leaves}
As mentioned in Section~\ref{ss:tropical} above, tropical geometers  call  metric graphs with leaves {\em tropical curves}, and the moduli space $\Mns$ the {\em tropical moduli space}.  From their point of view it is natural to compactify $\Mns$ by allowing subgraphs to collapse to points and recording the rank of the subgraph as an integer at the resulting vertex.  In our language, the compactification $\Mns^*$ is  the quotient of the simplicial closure $\Ons^*$ by $\G_{n,s}$.  A technical point here is that if we are only interested in the interior $\Mns$  we may retract all the separating edges to get a smaller space, since that does not change the homotopy type.  To study   $\Mns^*$, on the other hand, we may decide to keep graphs with separating edges; that is the compactification studied by tropical geometers. The space $\Mns^*$  can also be identified with the quotient of the curve complex of a surface of genus $n$ with $s$ punctures by the action of the mapping class group of the surface.  

Chan \cite{Chan} and  Chan, Galatius and Payne \cite{CGP} have studied the homology of $\Mns^*$ and shown in particular that it vanishes below dimension $s-3$.   Since $\Mns$ embeds in $\Mns^*$ one may wonder what happens under this embedding to the unstable homology classes we have found.     The answer is that any class in the image of assembly maps dies, because it is composed of classes supported on systems of subgraphs.  In $\Mns^*$ one may collapse any or all of the subgraphs in a system without collapsing the whole graph, thereby coning off (and killing) the entire assembled class.  This observation gives some (admittedly weak) credence to the conjecture made in \cite{CHKV} that all of the homology of $\Mns$ below dimension $2n+s-3$ is in the image of assembly maps.

\section{The geometry of Outer space - the Lipschitz metric}\label{sec:Lipschitz}
One way to measure  the difference between two marked metric graphs $(g,G)$ and $(g',G')$  in $\On$ is to find the map with minimal possible   Lipschitz constant among all homotopy equivalences $f\colon G\to G'$ with $f\circ g$ homotopic to $g'$.  This idea gives rise to the asymmetric {\em Lipschitz metric} on $\On$, whose basic properties were first detailed in \cite{FrMa}, and which has been explored in depth by Algom-Kfir, Bestvina and others \cite{Alg, AlBe}.  One striking application of this metric theory is Bestvina's streamlined ``Bers-like" proof  of the theorem that there is a very nice representative $f\colon G\to G$, called a {\em train track}, for a fully irreducible automorphism  of $F_n$ \cite{Bes1}.  Here $F_n$ is identified with $\pi_1(G)$ using the marking $g$.    Introductions to this  topic are available elsewhere, see, e.g. \cite{Bes2} or \cite{V_CE}.

\section{What's at infinity?}
A powerful technique in geometric group theory is to extend the action of a group on a space to an action on some nice compactification of that space, then study the action on the compactified space.    To study this action it is of course necessary to understand the points that were added, i.e. what is ``at infinity."   We start this section with a compactification of Outer space introduced at the same time as Outer space itself.
 
\subsection{Projective length functions}\label{ss:projective}
  Given a marked graph $(g,G)$ in $\On$ and a  word $w\in F_n= \pi_1(R_n)$, one can measure the length of the shortest loop in the homotopy class $g(w)$.  Since there is no basepoint to preserve, this depends only on the conjugacy class of $w$.  These lengths determine $(g,G)$ uniquely, so the entire space $\On$ embeds into   $(\R_{>0})^{\mathcal C}$ where $\mathcal C$ is the set of conjugacy classes in $F_n$.  Recall that we are only considering graphs with total volume one, so we can also think of this as a point in the projective space $\rm{P}\R^{\mathcal C}$.   Culler and Morgan showed that the closure of the image of $\On$ in  $\rm{P}\R^{\mathcal C}$ is compact \cite{CuMo}, and various authors  identified the points in the closure in terms of actions on trees \cite{CoLu, BeFe1, GaLe}. The actions which appear in the closure are not necessarily free and the trees are not necessarily simplicial; they may have dense branching.  The correct notion here is that of an $\R$-tree, and the actions of $F_n$ on $\R$-trees which appear in the closure are exactly the {\em very small} actions, where by definition an action is very small if  for each nontrivial $g\in F_n$ the fixed subtree $Fix(g)$ is isometric to a subset of $\R$ and is equal to $Fix(g^p)$ for all $p\geq 2$.

\subsection{Currents} \label{ss:currents}
Another compactification of $\On$ is formed by taking the closure of an embedding of $\On$ into the space of projectivized {\em geodesic currents} on $F_n$; these are $F_n$-invariant Borel measures on the set of unordered pairs of distinct points in the boundary of $F_n$. The definition of this embedding was motivated by work of Bonahon on Teichm\"uller space and first introduced into the free group context by I. Kapovich \cite{Kap}.  One important feature in the free group case is that  there is a natural {\em intersection pairing} between length functions and geodesic currents which extends to the Culler-Morgan boundary  and has found many applications.  One such application is an elementary proof that various simplicial complexes with an $\Outn$-action (such as the free factor complex and the simplicial closure of Outer space) have infinite diameter (recall from Section~\ref{ss:hyperbolic} that both of these complexes have since been shown to be Gromov hyperbolic.)  

\subsection{The  horoboundary and random walks}\label{ss:random} 
More recently another type of compactification, called the {\em horofunction compactification}, was studied by Horbez \cite{Hor1}.  Under mild conditions, which are satisfied by Outer space with its Lipschitz metric, a metric space $X$ can be embedded into the space $C(X)$ of continuous functions on $X$ as follows. Fix a basepoint $b\in X$, and to  each $z\in X$ associate the continuous function  $\psi_z(x)= d(z,x)-d(z,b)$; note that  the level sets of $\psi_z$ are spheres centered at $z$.    Given a geodesic ray $\zeta(t)$ in $X$ which leaves every compact set, there is a  {\em Busemann function} $B_\zeta$ defined by $$B_\zeta(x)= \lim_{t\to\infty} \psi_{\zeta(t)}(x);$$ the level sets of $B_\zeta$ are the horospheres centered at the limit point defined by $\zeta$.  The closure of the image of $X$ in $C(X)$ under the embedding $z\mapsto \psi_x$ is compact, the points in the boundary are called the {\em horoboundary}  and the points $B_\zeta$ in the horoboundary are called {\em Busemann points}.  

Outer space $\On$ embeds into the subspace of $\rm{P}\R^{\mathcal C}$  spanned by conjugacy classes of {\em primitive elements} of $F_n$; these are elements which can belong to a basis for $F_n$.  Horbez identified the horoboundary of Outer space with the closure of the image of $\On$ in this  subspace, then used this characterization in a new proof of the Tits alternative for  $\Outn$, which he generalized to a proof for the outer automorphism group of a free product \cite{Hor5}.  The idea is to use   the fact  that the complex $\FFn$ of free factors of $F_n$ is Gromov hyperbolic. By a theorem of Gromov there is a trichotomy in the possible behaviors of the action of a subgroup $H$ of $\Outn$ on $\FFn$:  either there are two loxodromic elements of $H$ which generate a free group, or there is a fixed point on the boundary of $\FFn$, or there is a bounded orbit.  The third case is the difficult one; one would like to conclude that in that case there is a fixed (conjugacy class) of free factors, and argue by induction.  But since $\FFn$ is not locally finite  one cannot reach this conclusion.  Horbez' solution is to prove that if there is no fixed free factor, then a random walk on the closure  of $\On$ produces a stationary measure, which he then uses to find a fixed point on the Gromov boundary of $\FFn$.  

Horbez's work on random walks culminates with a theorem echoing a classical theorem of Furstenburg for matrix groups \cite{Hor2}.  For $\Outn$ Horbez' theorem says that if a random walk is generated by a measure whose support generates all of $\Outn$, then  the length of almost every word $w\in F_n$ grows at a constant rate (depending only on $n$) under the random walk.    There is a more refined version of the theorem when the support of the measure generates a proper subgroup.  Finally,  a more detailed analysis of the random walk  gives rise to a ``central limit theorem" for $\Outn$, which describes the distribution of the average word length of primitive elements under iterated applications of random automorphisms (random with respect to a suitable measure)   \cite{Hor3}.

Horbez's investigations also result in  nice  a geometric description of the Poisson boundary of $\Outn,$ as the Gromov boundary of $\FFn$.  This is proved by following a point of Outer space under the action of a random sequence of automorphisms (with respect to some suitable measure on $\Outn$) and showing that it converges to a simplex in the boundary \cite{Hor3}.

\subsection{The Pacman compactification}  
Bartels, Lueck and Reich have introduced a geometric method of proving the Farrell-Jones conjecture for a group $G$.  This conjecture says that a certain map in $K$-theory is an isomorphism and it implies many other conjectures in manifold topology \cite{BLR}.  This method requires a compact contractible space with a proper $G$-action which is particularly well-behaved near the boundary.  
 
Both the Culler-Morgan compactification and the horofunction compactification of $\On$ have complicated local structure near the boundary, e.g. the boundary does not have a collar in the space (technically, the boundary is not a {\em Z-set}); this is a problem for the Bartels-Lueck-Reich method.   For $n=2$ the Culler-Morgan compactification is an absolute retract and it  is possible that one could get away with this weaker property,  but Bestvina and Horbez showed that   for $n\geq 4$ the Culler-Morgan compactification is not even an absolute neighborhood retract \cite{BeHo}.  They do this by showing  it is  not locally $4$-connected at a certain point which they describe explicitly.

In the same paper Bestvina and Horbez address these difficulties by defining a new compactification of $\On$. This a space of actions on $\R$-trees with additional structure, namely that arcs which have non-trivial stabilizers are assigned orientations.  They call this  the {\em Pacman compactification} because for $n=2$ adding orientations has the effect of slitting open the spikes at rational points in the boundary to form ``mouths" reminiscent of the classic Pacman video game.   The Pacman compactificaiton is an absolute retract  and the boundary is a Z-set, so it is a potential candidate for applying the Bartels-Lueck-Reich method.

\subsection{Bestvina-Feighn bordification}
In 2000 Bestvina and Feighn defined a {\em bordification} of $\On$  analogous to the Borel-Serre bordification of symmetric spaces for  non-compact semisimple algebraic groups defined over $\Q$ \cite{BeFe2}.  The Bestvina-Feighn construction follows that of Borel and Serre in spirit, embedding $\On$ into a larger space $\overline{\On}$ which is contractible (but not compact), and to which the action extends with compact quotient.  Their construction, and the  proof that the bordification is highly-connected at infinity, was intricate and left many questions unanswered, such as whether it was homeomorphic to a neighborhood of the spine of $\On$ and whether the building blocks were geometric cells.  

Grayson gave an alternate construction to that of Borel and Serre  in the case of the general linear group, finding a space homeomorphic to the bordification as a deformation retract of the symmetric space rather than an extension \cite{Grayson}. There is now an analogous picture for $\On$, namely a cocompact, equivariant deformation retract which forms a neighborhood of the spine, and which is highly-connected at infinity \cite{BuVo}.  The construction and proof that it is highly connected at infinity borrow heavily from Bestvina and Feighn's work, but give a clearer picture of the space and easier route to the proof.   

The construction of the bordification relies on the decomposition of  $\On$ into open simplices $\sigma(g,G)$, one for each isomorphism class of marked (combinatorial) graphs $(g,G)$ (see Section~\ref{ss:topology}).  Recall that faces of $\sigma(g,G)$ correspond to marked graphs obtained from $(g,G)$ by setting some of the edge lengths equal to zero, and, if the subgraph spanned by length 0 edges is not a forest,  we say the face is {\em at infinity}. In particular, 
if $G$ is a rose then all of the faces of $\sigma(g,G)$  are at infinity.  
 
 The bordification of $\On$ is a cell complex, all of whose vertices are contained in rose faces.  The intersection with each rose face is a well-known convex polytope called a {\em permutahedron}.   For example, in rank $3$ a rose face is a triangle, which we think of as equilateral.  After slicing off a neighborhood of each vertex and a smaller neighborhood of each edge, we are left with a hexagon, with one vertex for each permutation of the edges of the rose; this is the permutahedron for $n=3$.  
  For a general marked graph $(g,G)$ we have a permutahedron in each rose face of $\sigma(g,G)$; the convex hull in $\sigma(g,G)$ of all of their vertices is a closed cell $\Sigma(g,G)$.  The bordification $\overline{\On}$ is the union of the $\Sigma(g,G)$ for all marked graphs $(g,G)$.  To prove that this is highly-connected at infinity we apply a combinatorial Morse theory argument, where the Morse function is a function on the vertices which measures how many times $g(w)$ must pass over each edge of a rose $G$, for each conjugacy class $w\in F_n$.  The values of this Morse function lie in the ordered abelian semigroup $\mathbb N^\infty$, where $\mathbb N$ is the natural numbers and the ordering in lexicographical.

  \section{Outer connections}
  
  In this section we briefly touch on a few topics which connect $\On$ and moduli spaces of graphs with other topics in science and mathematics, with an indication of where to find further information.

  \subsection{Tropical moduli spaces}
  
  The connection of $\Ons$ to tropical moduli spaces was already mentioned in  {Section~\ref{ss:tropical}. Tropical geometers have been studying these moduli spaces, in particular establishing ``tropical" analogs of classic theorems from algebraic geometry about the moduli space of Riemann surfaces.  In particular they have studied the map from $\mathcal{MG}_{n,0}$ to the moduli space of flat tori,  which they call the {\em Torelli map}  (See \cite{Cap, BMV}).  
    
From our point of view it is more natural to define  this  map as an equivariant map from $\On$ to the symmetric space for ${\rm{SL}}(n)$, where it has traditionally been called the {\em Jacobian map} (see, e.g. the thesis of O. Baker \cite{Bak}).  Here the  definition  is very simple. To each marked graph $(g,G)$ we must associate a positive definite quadratic form on $\R^n$.   
To do this we use $g$ to identify $H_1(G;\R)$ with $H_1(R_n;\R)\cong \R^n$.   The first homology of $G$ is the kernel of the map from the $1$-chains to the $0$-chains, so can be thought of as a subspace of the vector space $\R^E$ with basis the edges $E$ of $G$.  There is a natural positive definite quadratic form on $\R^E$ in which the basis vectors are orthogonal and have length equal to the  length of the corresponding edge in $G$.  Restricting this form to $H_1(G;\R),$ which we have identified with $\R^n$, gives a positive definite quadratic form on $\R^n$, and hence a point in the usual symmetric space for ${\rm{SL}}(n)$.  
  
 There is at least one use for the Jacobian map which is closely connected to group theory.  The Jacobian descends to a map from the quotient space $\On/IA_n$ to the symmetric space, where   $IA_n$ is the kernel of the map from $\Outn$ to $\rm{GL}(n,\Z)$.   Since $IA_n$ is torsion-free, it acts freely on $\On$ and this quotient space is a genuine classifying space for $IA_n$.  The group $IA_n$ is quite mysterious, and this map can be used to investigate the structure of its classifying space, as Baker did in his thesis for the case $n=3$.   
 
 \subsection{Phylogenetic trees}  The space $\mathcal{MG}_{0,s}$ is the space of trees with $s$ labeled leaves. The leaves do not have lengths but the internal edges do, and we normalize so that the sum of the internal edge lengths is one. The cone $c\mathcal{MG}_{0,s}$ on this space is a non-positively curved (i.e. CAT(0)) metric space, and contains trees with all possible internal edge lengths.  
 
 Given a set of biological species which includes a common ancestor, we can label the leaves of a tree with the names of these species,  and thereby think of a point in $c\mathcal{MG}_{0,s}$  as a potential phylogenetic tree (i.e. a family tree, showing evolutionary relationships between species).  The interior edge lengths correspond to some measure of the distance between speciation events.    This is not an exact science:  measurements are only approximate, different measures may give different trees, and even if all parameters are agreed upon different runs of a computer program may spit out different trees because the the pairwise distances measured between species do not always correspond to distances in an actual tree.  One way to deal with this uncertainty is to generate a large number of possible trees and do a statistical analysis of the resulting  cloud of points in $c\mathcal{MG}_{0,s}$.    The fact that $c\mathcal{MG}_{0,s}$ is a  CAT(0) metric space allows one to  give several candidates for a meaningful average of a set of trees, and biologists are actively pursuing this idea (see, e.g. \cite{Owen}.)

 \subsection{Symplectic derivations of the free Lie algebra and number theory}
 
 In two seminal papers \cite{Kon1,Kon2} Kontsevich identified the cohomology of $\Outn$ with the homology of a certain infinite-dimensional symplectic Lie algebra $\ell_\infty$.  The Lie algebra $\ell_\infty$ is the direct limit of Lie algebras $\ell_n$ consisting of symplectic derivations of the free Lie algebra $L_n$ on $2n$ generators $\{p_1,\ldots,p_n,q_1,\ldots,q_n\}$.  Here a {\em symplectic derivation} is a linear map $D\colon L_n\to L_n$ that satisfies the Liebnitz rule $D[x,y]=[x, Dy] + [Dx, y]$ and vanishes on the  element $\omega_n=\sum [p_i,q_i]$.

  The same algebra $\ell_n$ appears in other mathematical contexts, including work of Morita and others on the Johnson filtration of the mapping class group of a surface \cite{Mor}, work of Garoufalidas and Levine on finite-type invariants of 3-manifolds \cite{GaLe}, and work of Berglund and Madsen on the rational homotopy theory of automorphisms of highly connected manifolds \cite{BeMa}.

 In fact Kontsevich described three diffeent symplectic Lie algebras $c_\infty, a_\infty$ and $\ell_{\infty}$, which he described as living in commutative, associative and Lie ``worlds" respectively.  The essential features of these worlds are captured in the notion of a cyclic operad.  Given any cyclic operad one can one can construct a symplectic Lie algebra  and a graph complex which computes the cohomology of this Lie algebra \cite{CoVo1}.

   \subsection{Feynman integrals}
 
 Feynman invented a method of computing quantum-mechanical amplitudes for a given physical system  by expanding them in a ``perturbative series" whose terms, called {\em Feynman integrals}, are indexed by finite graphs with leaves, known as {\em Feynman diagrams}.  The leaves of a Feynman diagram are labeled by momenta, which must sum to zero, and the edges by other parameters relating to the system, such as masses of particles and Schwinger normal times.  
 
There is a set of rules called {\em Cuttkosky rules} for calculating Feynman integrals in terms of Feynman diagrams over related graphs.  It turns out that the cubical structure of the spine of Outer space gives a natural way of organizing these related graphs, as detailed by Bloch and Kreimer in \cite{BlKr},  and it is  expected that the perspective of Outer space can further contribute to understanding these integrals.

\section{Previous survey articles}

There are a number of existing survey articles on aspects of topics discussed in the current survey. These may be useful for people who want to know more about a particular topic before they delve into the original sources.   These include  \cite{Bes2, Bes3,BrVo2, V_Hai,V_ICM,V_WI,V_CE, V_SA}.


\frenchspacing
\begin{thebibliography}{7}

% \bibitem{ }
% Author,
% Title.
%  \textit{Journal}  \textbf{volume} (date), pp--pp.

 
 
  \bibitem{Alg}% Yael's thesis?
Y. Algom-Kfir,
Strongly contracting geodesics in outer space.
\textit{Geom. Topol.}  \textbf{15}, no. 4, (2011), 2181--2233.

 \bibitem{AlBe}% some paper of  Yael and Mladen about the Lipschitz metric?
Y. Algom-Kfir and M. Bestvina,
Asymmetry of outer space. 
  \textit{Geom. Dedicata} \textbf{156 } (2012), 81--92.
 
  
 
 
\bibitem{Bak}% owen Baker Jacobian dissertation
O. Baker, 
The Jacobian map on Outer space.
Thesis, Cornell University (2011). 
Available at http://hdl.handle.net/1813/30769
 
 \bibitem{BLR}% Method for proving Farrell Jones
 A. Bartels, W. Lueck, and H. Reich, 
 The K-theoretic Farrell-Jones conjecture for hyperbolic groups.
  \textit{Invent. Math.}  \textbf{172} (2008), 29--70.
 
 \bibitem{Bar}% Bartholdi Out(F_7)
L. Bartholdi,
The rational homology of the outer automorphism group of $F_7$. 
\textit{New York J. Math.} \textbf{22} (2016), 191-197.
 
 
\bibitem{BeMa}
A. Berglund and I. Madsen, 
Rational homotopy theory of automorphisms of highly connected manifolds.
arXiv:1401.4096. 

 \bibitem{Bes1}% Bers-like proof of train track theorem
M. Bestvina,
A Bers-like proof of the existence of train tracks for free group
   automorphisms.
\textit{Fund. Math.} \textbf{214}  (2011),  no. 1, 1--12. 

 \bibitem{Bes2} % ICM talk
M. Bestvina,
The topology of $\Outn$. 
\textit{Proceedings of the International Congress of Mathematicians}, Vol. II (Beijing, 2002), 37--384, 
Higher Ed. Press, Beijing, 2002.
 
  
  \bibitem{Bes3} % Park city notes
M. Bestvina,
Geometry of outer space.
\textit{Geometric group theory},
173--206, \textit{IAS/Park City Math. Ser.} \textbf{21}, Amer. Math. Soc., Providence, RI,  2014.
 
\bibitem{BBF} M. Bestvina, K. Bromberg, and K. Fujiwara.
 Constructing group actions on quasi-trees and applications to mapping class groups.
\textit{Publ. Math. Inst. Hautes ƒtudes Sci.} \textbf{122} (2015), 1--64. 

       
 \bibitem{BeFe}% Hyperbolicity of free factor complex
M. Bestvina and M. Feighn,
Hyperbolicity of the complex of free factors. 
\textit{Adv. Math.} \textbf{256} (2014), 104--155

 \bibitem{BeFeErr}% Hyperbolicity of free factor complex
M. Bestvina and M. Feighn,
Corrigendum to "Hyperbolicity of the complex of free factors'' [\textit{Adv. Math.} \textbf{256} (2014) 104--155].
\textit{Adv. Math.} \textbf{259} (2014), 843.
 
 \bibitem{BeFe1}%Outer limits
M. Bestvina and M. Feighn, Outer limits.
(Unpublished.) 
Available in preprint form at 
http://andromeda.rutgers.edu/\textasciitilde feighn/papers/outer.pdf
 
  \bibitem{BeFe2}%Bordification
M. Bestvina and M. Feighn,
The topology at infinity of $\Outn$. 
\textit{Invent. Math.} \textbf{140} (2000), no. 3, 651--692.
 

 
\bibitem{BeHo}
M. Bestvina and C. Horbez,
 A compactification of outer space which is an absolute retract. 
arXiv:1512.02893.

 
 \bibitem{BiHoVo}
L. Billera, S. Holmes and K. Vogtmann
Geometry of the space of phylogenetic trees. 
 \textit{Adv. in Appl. Math.} \textbf{27} (2001), no. 4, 733--767. 

 \bibitem{BlKr}%Cuttkosky rules
S. Bloch and D. Kreimer,
 Cutkosky Rules and Outer Space.
     arXiv:1512.01705.

\bibitem{BMV} 
S. Brannetti, M. Melo and F. Viviani,
On the tropical Torelli map.
\textit{Adv. Math.} \textbf{226} (2011),
no. 3, 2546--586.


\bibitem{BrVo}
 M. R. Bridson and K. Vogtmann, The Dehn functions of $\Outn$ and $\Autn$. 
 \textit{Ann. Inst. Fourier(Grenoble)} \textbf{62} (2012), no. 5, 1811--1817.
 
\bibitem{BrVo2}
M. R.  Bridson and K. Vogtmann, 
Automorphism groups of free groups, surface groups and free abelian groups. 
\textit{Problems on mapping class groups and related topics}, 301--316, 
\textit{Proc. Sympos. Pure Math.} \textbf{74},
Amer. Math. Soc., Providence, RI, 2006.



\bibitem{BuVo}
K.-U. Bux and K. Vogtmann,
in preparation.

\bibitem{Cap}% Caporoso?? etc on tropical curves.get reference from Chan? Then alphabetize!
L. Caporaso, 
Algebraic and tropical curves: comparing their moduli spaces,
\textit{Handbook of Moduli, Vol. I},  
 \textit{Adv. Lect. Math.}  \textbf{24}  Int. Press, Somerville, MA, 2013, 119-160.
 
 
 
 \bibitem{Chan}
M. Chan,
Topology of the tropical moduli spaces $M_{2,n}$.
 arXiv:1507.03878 .
 
\bibitem{CGP}
M. Chan, S. Galatius and S. Payne
The tropicalization of the moduli space of curves II: Topology and applications.
 arXiv:1604.03176.

\bibitem{ChVo} 
R. Charney and K. Vogtmann,
Outer space for untwisted automorphism groups of RAAGs.
To appear in \textit{Geom. and Topol.}

\bibitem{CoLu}%Very small actions
M. Cohen and M. Lustig,
Very small group actions on R-trees and Dehn twist automorphisms. 
\textit{Topology} \textbf{34} (1995), no. 3, 575--617. 


\bibitem{CoVo1}
J. Conant and K. Vogtmann,
On a theorem of Kontsevich.
 \textit{Algebr. Geom. Topol.}  \textbf{3} (2003), 1167--1224.
 
 
\bibitem{CoVo2}
J. Conant and K. Vogtmann,
Morita classes in the homology of automorphism groups of free groups. 
 \textit{Geom. Topol.}  \textbf{8} (2004), 1471--1499 

\bibitem{CKV}
J. Conant, M. Kassabov and K. Vogtmann,
Hairy graphs and the unstable homology of $\rm{Mod}(g,s)$, $\Outn$ and $\Autn$.
 \textit{J. Topol.}  \textbf{6} (2013), no. 1, 119--153.

\bibitem{CHKV}
J. Conant, A. Hatcher, M. Kassabov and K. Vogtmann,
Assembling homology classes in automorphism groups of free groups,
to appear in \textit{Comment. Math. Helv.}   

\bibitem{CuMo}
M. Culler and J. Morgan,
Group actions on R-trees. 
\textit{Proc. London Math. Soc.} (3) \textbf{55} (1987), no. 3, 571--604.
 
 
\bibitem{CuVo}
M. Culler  and  K. Vogtmann,
Moduli of graphs and automorphisms of free groups.
\textit{Invent. Math}  \textbf{84} (1986), 91--119.

\bibitem{FrMa} %  Lipschitz metric
S. Francaviglia and A. Martino,
Metric properties of outer space.
\textit{Publ. Mat.}  \textbf{55}  (2011),  no. 2, 433--473.
  
 
 \bibitem{GaLe}%Boundary of OS
D. Gaboriau and G. Levitt
The rank of actions on R-trees. 
\textit{Ann. Sci. ƒcole Norm. Sup.} (4) \textbf{28} (1995), no. 5, 549--570.


 \bibitem{Gal}%  Stable homology of Out(F_n)
S. Galatius
 Stable homology of automorphism groups of free groups.
\textit{Ann. of Math.} (2)  \textbf{173}  (2011),  no. 2, 705--768.
 
\bibitem{GaLe}% Finite-type invariants...look at CKV for reference
 S.  Garoufalidis  and  J.  Levine,
   Tree-level  invariants  of  $3$-manifolds,  Massey  products  and
the Johnson homomorphism. 
\textit{Graphs and patterns in mathematics and theoretical physics} 173--203, 
 \textit{Proc. Sympos. Pure Math.} \textbf{73}, Amer. Math. Soc., Providence, RI, (2005).
 
    
  \bibitem{Ger}% Not published.   
F. Gerlits,
Invariants in Chain Complexes of Graphs.
Thesis, Cornell University  2002.   

 \bibitem{Gra} % Not published.
J. Gray, 
On the homology of automorphism groups of free groups.
Thesis, University of Tennessee, 2011. Available at 
http://trace.tennessee.edu/utk-grad/  

\bibitem{Grayson}
D. Grayson,
Reduction theory using semistability.
\textit{Comment. Math. Helv.}  \textbf{59}  (1984),  no. 4, 600--634.


 \bibitem{GiLe1}% paper on automorphisms of free products
V. Guirardel and G. Levitt,
The outer space of a free product.
\textit{Proc. Lond. Math. Soc.} (3)  \textbf{94}  (2007),  no. 3, 695--714.

\bibitem{GiLe2}% paper on deformation spaces
V. Guirardel and G. Levitt,
Deformation spaces of trees.
\textit{Groups Geom. Dyn.}  \textbf{1}  (2007),  no. 2, 135--181.

\bibitem{HaMo1}% Dehn function for Out.
M. Handel, L. Mosher,
Lipschitz retraction and distortion for subgroups of $\Outn$.
\textit{Geom. Topol.} \textbf{17} (2013), no. 3, 1535--1579.


\bibitem{HaMo2}% Hyperbolicity of free splitting complex
M. Handel, L. Mosher,
The free splitting complex of a free group, I: hyperbolicity.
\textit{Geom. Topol.}  \textbf{17}  (2013),  no. 3, 1581--1672.
 

\bibitem{Hat}% Homology stability paper (sphere complexes}
A. Hatcher,
Homological stability for automorphism groups of free groups.
\textit{Comment. Math. Helv.} \textbf{70}  (1995),  no. 1, 39--62.
 
 \bibitem{HaVo0}%  Rational homology (cube complex structure)
 A. Hatcher and K. Vogtmann
 Rational homology of $\Autn$.
\textit{ Math. Res. Lett.} \textbf{5} (1998), no. 6, 759--780.
 
 \bibitem{HaVo1}% Dehn function paper
 A. Hatcher and K. Vogtmann,
Isoperimetric inequalities for automorphism groups of free groups. 
 \textit{Pacific J. Math.}  \textbf{173} (1996), no. 2, 425--441.  
 
  \bibitem{HaVo2}%Homology stability
A. Hatcher and K. Vogtmann,
Homology stability for outer automorphism groups of free groups.
 \textit{Algebr. Geom. Topol.}  \textbf{4} (2004), 1253--1272.
 


  \bibitem{HVW}%Out   stability erratum
A. Hatcher, K. Vogtmann and N. Wahl,
Erratum to: "Homology stability for outer automorphism groups of free groups.     
 \textit{Algebr. Geom. Topol.}  \textbf{6} (2006), 573--579.
 
  \bibitem{HaWa}% Paper on homology stability for 3 manifolds (plus erratum?)
A. Hatcher and N. Wahl,
 Stabilization for mapping class groups of 3-manifolds.
\textit{Duke Math. J.}  \textbf{155}  (2010),  no. 2, 205--269.


  \bibitem{HiHo}% Hyperbolicity of sphere complex
A. Hilion and C. Horbez,
The hyperbolicity of the sphere complex via surgery paths,  
 to appear in J. Reine Angew. Math., arXjiv:1210.6183.
 
 \bibitem{Hor1}% Horoboundary
 C. Horbez,
The horoboundary of outer space, and growth under random automorphisms. 
 arXiv:1407.3608.
 
 \bibitem{Hor2}% Central limit
 C. Horbez,
Central limit theorems for mapping class groups and $\Outn$.
 to appear in Ann. Scient. Ec. Norm. Sup.,  arXiv:1506.07244.

 \bibitem{Hor3}% Poisson
 C. Horbez,
 The Poisson boundary of $\Outn$.
 \textit{Duke Math. J.} \textbf{165} (2016), no. 2, 341--369.
 
 \bibitem{Hor5}% Tits alternative
 C. Horbez,
The Tits alternative for the automorphism group of a free product.
arXiv:1408.0546.

\bibitem{Kap}% embedding  of OS into currents
I. Kapovich,
Currents on free groups.
Topological and asymptotic aspects of group theory,
149--176, 
\textit{Contemp. Math.} \textbf{394}, Amer. Math. Soc., Providence, RI, 2006.

\bibitem{Kon1}
M. Kontsevich,
Feynman diagrams and low-dimensional topology. 
\textit{First European Congress of Mathematics, Vol. II  (Paris, 1992)}, 97--121,
 \textit{Progr. Math.} \textbf{120}, BirkhŠuser, Basel, 1994.

\bibitem{Kon2}%The other Kontsevich paper...
M. Kontsevich, 
Formal (non)commutative symplectic geometry. 
\textit{The Gelfand Mathematical Seminars, 1990--1992},
173--187, BirkhŠuser Boston, Boston, MA, 1993.
 
 \bibitem{Lau}
F. Laudenbach,
 Sur les 2-sphres d'une variŽtŽ de dimension 3. 
 \textit{Ann. of Math.} (2) \textbf{97} (1973), 57--81.
 
 

 
  \bibitem{Mor}% Survey and prospects paper
S. Morita,  
Structure of the mapping class groups of surfaces: a survey and a prospect.
\textit{Proceedings of the Kirbyfest (Berkeley, CA, 1998)},
349--406 (electronic), \textit{Geom. Topol. Monogr.} 2, Geom. Topol. Publ., Coventry, 1999.

\bibitem{MSS}% Euler characteristic n=11.
S. Morita, Sakasai, Suzuki, 
Integral Euler characteristic of $\Out(F_{11})$. 
\textit{Exp. Math.} \textbf{24} (2015), no. 1, 93--97. 

\bibitem{MSS2} % Euler characteristic n<11.
 S. Morita, T. Sakasai, M. Suzuki, 
Computations in formal symplectic geometry and characteristic classes of moduli spaces.
\textit{Quantum Topol.} \textbf{6} (2015), no. 1, 139--182.
 
 \bibitem{Osi}%Osin survey about acylindrical groups?
 D. Osin,
 Acylindrically hyperbolic groups. 
 \textit{Trans. Amer. Math. Soc.} \textbf{368} (2016), no. 2, 851--888.

 \bibitem{Owen}% Owens work on phylogenetic trees
E. Miller, M. Owen, and J. S. Provan,
Polyhedral computational geometry for averaging metric phylogenetic trees. 
\textit{Adv. in Appl. Math.} \textbf{68} (2015), 51--91. 



\bibitem{Sko}
R. Skora,
Deformations of length functions in groups.
 (Unpublished.)  Available at http://www.math.unicaen.fr/~levitt/unpublished.html
 
 \bibitem{SmVo}
 J. Smillie and K. Vogtmann,  
 A generating function for the Euler characteristic of $\Outn$. 
 Proceedings of the Northwestern conference on cohomology of groups (Evanston, Ill., 1985). 
 \textit{J. Pure Appl. Algebra} \textbf{44} (1987), no. 1-3, 329--348.

\bibitem{V_Hai}% Haifa paper
K. Vogtmann,
Automorphisms of free groups and outer space.
\textit{Proceedings of the Conference on Geometric and Combinatorial Group Theory, Part I} (Haifa, 2000).
 \textit{Geom. Dedicata}  \textbf{94} (2002), 1--31.


\bibitem{V_ICM}% ICM paper
K. Vogtmann,
The cohomology of automorphism groups of free groups.
\textit{Proceedings of the International Congress of Mathematicians. Vol. II}  (2006), 1101--1117, 
 Eur. Math. Soc., Z\"urich.


\bibitem{V_WI}% What is
  K. Vogtmann,
  What is $\ldots$ Outer space? 
 \textit{Notices Amer. Math. Soc.}  \textbf{55} (2008) no. 7, 784--786. 



\bibitem{V_CE}% Current Events paper
 K. Vogtmann,
On the geometry of Outer space.
 \textit{Bull. Amer. Math. Soc. (N.S.)}  \textbf{52}, no. 1 (2015), 27--46.
 
 \bibitem{V_SA}% St. Andrews
K. Vogtmann,
$\rm{GL}(n, \Z)$, $\Outn$ and everything in between: automorphism groups of RAAGs, 
\textit{Groups St Andrews 2013}, 105--127, \textit{London Math. Soc. Lecture Note Ser.}  \textbf{422}: 
 Cambridge University Press, 2015. 

 
 \bibitem{V_Rep} % Reprise paper
 K. Vogtmann,
 Contractibility of Outer space: reprise.
Contractibility of Outer space:reprise, to appear in Math Society of Japan Summer Institute 2014 Proceedings,  arXiv:1505.02610.

\bibitem{Whi} %Tad White:  contractibility of fixed point set
 T. White,
Fixed points of finite groups of free group automorphisms. 
\textit{Proc. Amer. Math. Soc.} \textbf{118} (1993), no. 3, 681--688.
 
 
\end{thebibliography}
\end{document}